\documentclass[a4paper, 12pt]{article}

\usepackage[unicode]{hyperref}
\usepackage{amsmath}
\usepackage{amssymb}
\usepackage[strict,style=british]{csquotes}
\usepackage{graphics}
\usepackage{booktabs}
\usepackage{geometry}
\usepackage{authblk}
\usepackage[sort]{cite}
\usepackage{tablefootnote}
\title{Model reduction in Smoluchowski-type equations}
\author[1,2,3]{I.~V.~Timokhin\thanks{m@ivan.timokhin.name}}
\author[4,2]{S.~A.~Matveev}
\author[2,1,3]{E.~E.~Tyrtyshinikov}
\author[1]{A.~P.~Smirnov}
\date{}

\affil[1]{Faculty of Computational Mathematics and Cybernetics, Lomonosov Moscow State University, Russia}
\affil[2]{Marchuk Institute of Numerical Mathematics of Russian Academy of Sciences, Moscow, Russia}
\affil[3]{Moscow Center for Fundamental and Applied Mathematics, Moscow, Russia}
\affil[4]{Skolkovo Institute of Science and Technology, Moscow, Russia}

\begin{document}

\maketitle{}

\begin{abstract}
In this paper we utilize the Proper Orthogonal Decomposition (POD) method for model 
order reduction in application to Smoluchowski aggregation equations with source
and sink terms. 
In particular, we show in practice that there exists a low-dimensional space
allowing to approximate  the solutions of aggregation equations. 
We also demonstrate that 
 it is possible to model the aggregation process with the complexity depending only on dimension
 of such a space but not on the original problem size. In addition, we propose a method for
 reconstruction of the necessary space without solving of the full evolutionary problem, which can lead to
 significant acceleration of computations, examples of which are also presented.


\end{abstract}

\section{Introduction}\label{sec:intro}

A classical model of aggregation kinetics is based on the
Smoluchowski equations, dating back to the original work by Marian
von~Smoluchowski~\cite{smoluchowski1916drei}.  In the original form,
these equations describe an evolution of a spatially uniform system of agglomerates
of different sizes, via an infinite system of ordinary differential equations for
the concentrations $n_k$ of particles of size $k$ each.  The original
formulation has been later amended by Hans Muller~\cite{muller1928allgemeinen}
to model continuous particle size distribution or additional phenomena, such as
particle fragmentation~\cite{blatz1945note} and others.

The range of phenomena modelled via Smoluchowski kinetic equations has also
expanded over time, from molecular scales~\cite{blatz1945note, Privman1999, boje2019hybrid, boje2020study} to
astronomical~\cite{PNAS, esposito2012predator, esposito2008moonlets}. More detailed 
information about possible applications of aggregation-based models can be found in extensive 
reviews~\cite{leyvraz2003scaling, semeniuk2020current} and references therein.

Whether an original discrete system is used, or a discretization of the
continuous, one still has to deal with a rather large systems of
nonlinear differential equations, especially if particle
masses differ by several orders of magnitude.  Because of
that, accurate numerical simulation of these systems is quite
challenging.  While there has been some recent progress on this front,
bringing complexity for some classes of coagulation kernels down to
almost linear~\cite{chaudhury2014computationally, timokhin2019tensorisation} it is still
insufficient for some of the larger systems arising in practice.

In this paper, we will attack this problem using the ideas
of model reduction via Proper Orthogonal
Decomposition, as outlined in~\cite{modelpod}. Specifically, we
are interested in the method of snapshots,
introduced in~\cite{turbulence}.  The main idea of the method is to
construct a low-dimensional vector space containing the solution or its approximation
by examining its snapshots at different time moments.  The end goal here
is to create an opportunity to describe and approximate the solution using
significantly fewer parameters than the full dimensionality of the system.

In this paper we demonstrate that
\begin{itemize}
\item a low-dimensional space in which the solution can be approximated
 with reasonable accuracy \emph{exists};
\item once such a space is found, it is possible to model the system within the
  complexity depending only its dimension but not on the original problem size;
\item at least for some cases, it is possible to find the necessary space
  without constructing the solution of the full original problem.
\end{itemize}

Even though the results presented here do not seem to be immediately applicable
 for complex industrial applications, we believe that we suggest a
 novel concept for solving the aggregation-fragmentation equations leading to a
 fruitful and challenging avenue of further research. In some sense, our approach gives an alternative
 view at developing deep learning-based methods~\cite{shalova2020deep} for non-linear 
 time-dependent problems with attractors and cycles. In contrast to~\cite{shalova2020deep} 
 we deal with much larger systems of ODEs (tens of thousands in our work instead of 
 dozens or hundreds).
 
 The rest of the paper is organized as following: in Section 2 we discuss the target set 
 of kinetic equations and recall the necessary facts about properties of the solution 
 and the coefficients. In Section 3 we introduce a numerical method allowing to solve
 the target equations in approximate form using the reduction basis. The next Section 4
 is devoted to algorithm allowing to construct such a basis via Proper Orthogonal 
 Decomposition (namely, the method of snapshots). In Section 5 we demonstrate
 the results of numerical experiments and validation of the proposed methodology. In
 our experiments, we demonstrate the existence of the required low-dimensional reduced
 basis allowing one to accelerate the computations of numerical solutions of aggregation
 equations. In this Section, we also discuss the drawbacks of the proposed approach and 
 further accumulate our findings in the conclusions of Section 6.

\section{Problem setting}\label{sec:problem}

In our work, we consider the model similar to one originally posed
in~\cite{smoluchowski1916drei}, with the addition of a constant
source of particles~\cite{hayakawa1987irreversible, ball2012collective}:
\begin{equation}
  \label{eq:smoluchowski_inf}
  \frac{dn_k}{dt}
  = J_k
  + \frac{1}{2}\sum_{i + j = k} C_{i\,j} n_i n_j
  - n_k \sum_{j = 1}^{\infty} C_{j\,k} n_j,
  \qquad k = \overline{1, \infty}
\end{equation}

In this system,
\begin{description}
\item[$n_k$] stands for the concentration of particles of mass $k$;
\item[$C_{i\,j}$] is a coagulation kernel, characterising the frequencies
  of collisions between particles of size $i$ and $j$;
\item[$J_k$] is a uniform source of particles of size $k$.
\end{description}

We additionally put some physically relevant constraints on these variables:
\begin{description}
\item[$n_k \geq 0$] (there cannot be a negative concentrations of any kind of
  particles in the system);
\item[$C_{i\,j} = C_{j\,i} \geq 0$] (the coagulation kernel is symmetric
  and non-negative);
\item[$J_k \geq 0$] (this term corresponds to the \emph{source} of new particles).
\end{description}

For modelling purposes, we truncate~(\ref{eq:smoluchowski_inf}) to get
a finite system; this is equivalent to postulating an immediate
removal of large particles from the system (see
e.g.~\cite{ball2012collective}):
\begin{equation}
  \label{eq:smoluchowski}
  \frac{dn_k}{dt}
  = J_k
  + \frac{1}{2}\sum_{i + j = k} C_{i\,j} n_i n_j
  - n_k \sum_{j = 1}^{N} C_{j\,k} n_j,
  \qquad k = \overline{1, N}
\end{equation}

Given a sufficiently large $N$, system~(\ref{eq:smoluchowski})
approximates~(\ref{eq:smoluchowski_inf}) with reasonable accuracy
either in steady-state form~\cite{matveev_anderson_2017} or
quasi-steady-state~\cite{timokhin2019newton}.  For some cases of
kernel coefficients with finite $N$ a steady collective
oscillatory solutions of aggregation
equations~\cite{ball2012collective} exists, which cannot be expected for the
pure infinite aggregation system with source but no sink.
However, the required value of $N$ in practice can still be fairly
large, so our aim for the rest of the paper is to reduce the number of
parameters in~(\ref{eq:smoluchowski}).

\section{Model reduction}\label{sec:reduction}

In order to reduce the number of variables in the
system~(\ref{eq:smoluchowski}), we employ the model reduction concept
via Proper Orthogonal Decomposition (POD)~\cite{modelpod}.  The \emph{output} of
the method is an orthonormal basis of a low-dimensional subspace containing 
the solution allowing  at least to construct its approximation.  For now,
let us assume that we have already found the basis, and see how it can help
to work with the aggregation equations~(\ref{eq:smoluchowski}).

Let us start by rewriting~(\ref{eq:smoluchowski}) in a more
general form.  Namely, we start by introducing a tensor
$S \in \mathbb{R}^{N \times N \times N}$:
\begin{equation}
  \label{eq:smoluchowski_tensor}
  S_{i\,j\,k} = \frac{1}{2}\left(\delta_{i+j,k} - \delta_{i,k} - \delta_{j,k}\right) C_{i\,j},
\end{equation}
where $\delta_{i,j}$ is the Kronecker symbol.  Armed with this tensor,
we rewrite~(\ref{eq:smoluchowski}) as
\begin{equation}
  \label{eq:smoluchowski_tensor_eqs}
  \frac{dn_k}{dt} = J_k + \sum_{i,j = 1}^{N} S_{i\,j\,k} n_i n_j.
\end{equation}

Further, we assume the existence of an orthonormal basis, gathered as columns of a
matrix $V \in \mathbb{R}^{N \times R}$, such that
\begin{equation}
  \label{eq:basis_approx}
  \left \lVert n(t) - VV^T n(t) \right\rVert \ll \left \lVert n(t) \right \rVert,
\end{equation}
where $n(t)$ is the solution to~(\ref{eq:smoluchowski_tensor_eqs}).
We will hereafter abbreviate inequalities of this sort to
$n(t) \approx VV^T n(t)$.

Then we can introduce
\begin{equation}
  \label{eq:Vn}
  x(t) \equiv V^T n(t),\qquad x(t) \in \mathbb{R}^R,
\end{equation}
so that the equation~(\ref{eq:basis_approx}) turns into $n(t) \approx Vx(t)$. 
Substituting it into the equation~(\ref{eq:smoluchowski_tensor}), we get
\begin{equation}
  \label{eq:smol_V_approx}
  \frac{d}{dt} \sum_{\alpha=1}^R V_{k\,\alpha} x_\alpha(t)
  \approx J_k
  + \sum_{i,j=1}^N S_{i\,j\,k}
  \times
  \left(
    \sum_{\beta=1}^R V_{i\,\beta} x_\beta(t)
  \right)
  \times
  \left(
    \sum_{\gamma=1}^R V_{j\,\gamma} x_\gamma(t)
  \right).
\end{equation}

Multiplying this last system by $V^T$ and rearranging the sums a bit
we arrive at a reduced form of the original system (note that doing so
does not increase the second-norm absolute error, although it may well
increase the relative one):
\begin{equation}
  \label{eq:smol_reduced_full}
  \frac{d}{dt} x_\alpha
  \approx \sum_{k=1}^N V_{k\,\alpha} J_k
  + \sum_{\beta,\gamma=1}^R
  \left(
    \sum_{i,j,k=1}^N S_{i\,j\,k}
    V_{i\,\alpha}
    V_{j\,\beta}
    V_{k\,\gamma}
  \right)
  x_\beta
  x_\gamma,
\end{equation}
or introducing some extra notation
\begin{align}
  \label{eq:Jred}
  \tilde{J}_\alpha &= \sum_{k=1}^N V_{k\,\alpha} J_k, \\
  \label{eq:Sred}
  \tilde{S}_{\alpha\,\beta\,\gamma} &=
                                      \sum_{i,j,k=1}^N S_{i\,j\,k}
    V_{i\,\alpha}
    V_{j\,\beta}
    V_{k\,\gamma},
\end{align}
we rewrite it as
\begin{equation}
  \label{eq:smol_reduced_approx}
  \frac{d}{dt} x_\alpha \approx \tilde{J}_\alpha + \sum_{\beta,\gamma = 1}^R \tilde{S}_{\alpha\,\beta\,\gamma} x_\beta x_\gamma
\end{equation}

Finally, instead of defining $x$ via $n$, we can
recast~(\ref{eq:smol_reduced_approx}) as a system of ordinary
differential equations for a new variable $\tilde{x}$, that approximates
$x$:
\begin{align}
  \label{eq:smol_reduced}
  \frac{d}{dt} \tilde{x}_\alpha &= \tilde{J}_\alpha + \sum_{\beta,\gamma = 1}^R \tilde{S}_{\alpha\,\beta\,\gamma} \tilde{x}_\beta \tilde{x}_\gamma,\qquad \alpha = 1,2,\dots{},R \\
  \tilde{x}(0) &= x(0) = V^T n(0).
\end{align}

The important thing to notice here is that evaluation of the right-hand
part of~(\ref{eq:smol_reduced}) only takes $O(R^3)$ operations. Hence,
we reach our initial goal of completely decoupling the dimensionality of
the reduced system from $N$.
If $ R \ll N$ it may lead to a significant speedup of computations.

The reduced solution $\tilde{x}$ can then be used to reconstruct an
approximation to the full solution by further approximating the original equation~(\ref{eq:basis_approx}):
\begin{equation}
  n(t) \approx \tilde{n}(t) = V \tilde{x}(t).
\end{equation}

\section{Constructing a basis}\label{sec:basis}

In this section we describe a method used to construct the basis $V$ which 
we have been using in the previous section. To fullfill this aim we use the 
snapshot method from~\cite{modelpod}.

In this method, the basis $V$ is constructed via the snapshots of the
original solution at some fixed moments in time $t_k$, for $k = 1,\dots{},m$.
The exact method of basis construction may vary in technical details in different
publications about its applications, but the specific method from~\cite{modelpod} 
ends up being equivalent to taking leading left singular vectors of the $N \times m$ 
matrix of snapshots.  Specifically, in our case, $V$ is taken to be a matrix of senior left singular
vectors of a matrix composed of \enquote{snapshots} $n(t_k)$, with time moments $t_k$ uniformly spaced
across the interval of interest.  The number of singular vectors depends on the
specified approximation requirements; in practice, we use the same criteria as when combining bases (see below).

Unfortunately, we essentially need to \emph{know} the solution for 
construction a reduced basis which we are going to
to use to find of the approximation of the solution
.  To resolve this circularity,  we
split the initial time-interval into a number of \enquote{windows}, and use the snapshot
method to construct a basis for each of them in turn 
instead of finding just one basis for the entire time segment of our interest.

Specifically, let $\tau$ be some fixed time-window width, and assume we
have $\hat{V}_k$ such that
\begin{align}
  n(t) \approx \hat{V}_k \hat{V}_k^T n(t),\qquad \forall t \in [(k-1)\tau, k\tau].
\end{align}

Each of these can be constructed via the method of snapshots by
numerically solving of the full system~(\ref{eq:smoluchowski}) at  each
\enquote{window} in turn by use of any standard numerical method for ODE systems.

To combine them into a final, common basis, we introduce an auxiliary
operation $\oplus_\delta$ for any given $\delta > 0$: given two
matrices $A \in \mathbb{R}^{N \times r_1}$ and
$B \in \mathbb{R}^{N \times r_2}$, $A \oplus_\delta B$ is a
$N \times r_3$ matrix composed of the senior $r_3$ left singular
vectors of an $N \times (r_1 + r_2)$ matrix $C = (A \mid B)$ (that is,
a matrix composed of columns of $A$ and $B$~--- in principle, in any
order), where $r_3$ is chosen so that
$\sigma_{r_3} \geq \delta > \sigma_{r_3 + 1}$, where $\sigma_k$ are
singular values of $C$.

As a measure of the quality of our basis, we measure an error of
approximation of the \emph{next} window's basis by the ``current'' one, with some
small positive tolerance $\varepsilon \ll 1$.  In other words, our algorithm for
the basis construction can be formulated as following:
\begin{description}
\item[\textbf{Step 1}] Set $k \leftarrow 1$, $V_0 = 0 \in \mathbb{R}^{N \times 0}$.
\item[\textbf{Step 2}] Calculate $\hat{V}_k$ via the method of snapshots as an approximate reduction basis for the time span $[(k - 1)\tau, k\tau]$.
\item[\textbf{Step 3}] If
  $\left\lVert (I - V_{k-1} V_{k-1}^T) \hat{V}_k \right\rVert_2 \leq
  \varepsilon$, set $V = V_{k-1}$ and exit.
\item[\textbf{Step 4}] Otherwise, set $V_k = V_{k-1} \oplus_\delta \hat{V}_k$.
\item[\textbf{Step 5}] Set $k \leftarrow k + 1$ and repeat from step 2.
\end{description}

This algorithm is \enquote{greedy} in some sense: it tries to
approximate the entire solution by aggressively approximating each
subsequent time \enquote{window}. Hence, it probably may lead
to overestimation of the eventual basis dimensionality. In our concrete 
implementation we execute \textbf{Step 4} only if the approximation error at
\textbf{Step 3} is larger than an additional  auxiliary parameter 
$\varepsilon' > \varepsilon > 0$.

\section{Numerical experiments}\label{sec:numeric}

In this section we present the tests of the implementation of our method 
from previous sections. In our simulations we use a classical Browninan-type kernel
\begin{equation}
  \label{eq:kernel}
  C_{i\,j} = i^a j^{-a} + i^{-a} j^a.
\end{equation}
Even though problems with such kernel and its closest generalizations 
\[
C_{i\,j} = i^{\nu} j^{\mu} + i^{\mu} j^{\nu} + 2
\]
are rather well-studied by nowadays~\cite{hayakawa1987irreversible, krapivsky2012driven,
ball2012collective, matveev2020oscillating, krapivsky2018steady127691588}
the exact analytical solutions for time-dependent cases are still unknown especially for the cases with 
steady oscillations~\cite{ball2012collective, krapivsky2018steady127691588}. Moreover, researchers are still interested in the exploitation of such kernels 
 for practical modelling~\cite{slomka2020bursts} and theoretical analysis~\cite{pego2020temporal} as well.

In~\cite{matveev2015a-fast8543937}, a fast numerical method is
given for a Cauchy problem with this kernel, evaluating the
right-hand side of the equation~(\ref{eq:smoluchowski}) in just $O(N \log N)$
 operations and we want to out-perform this approach.
In~\cite{ball2012collective, krapivsky2018steady127691588},
steady collective oscillations in time of $n_k$  were detected for systems with 
this family of kernels with $a > 0.5$.

We have chosen this system specifically because the presence of the
cycles in the solution gives us a strong a~priori reason to expect that
our method works.  Namely, at least after the first iteration of the
cycle, any basis which adequately approximates the solution should
also approximate the further solution and can be used
at least to verify the cyclical behaviour.  However, we also note that due
to the use of two thresholds $\varepsilon' > \varepsilon > 0$ we
cannot \emph{guarantee} that the algorithm terminates, but, as
we soon see in our experiments, it frequently does.

At first, we demonstrate the principle feasibility and inner workings
of the algorithm.  For this purpose, we present a sequence of
experiments with the following set of model parameters:
\begin{align*}
  N &= 32768,&
  J_k &= \delta_{k\,1},\\
  \tau &= 2,&
  \varepsilon &= 10^{-13},\\
  \varepsilon' &= 10^{-10},&
  \delta &= 10^{-13}, \\
\end{align*}
where $\delta$ is for our \enquote{basis addition operator} $\oplus_\delta$ from Section~\ref{sec:basis}, with
$m = 65$ snapshots in each window for the snapshot method.  As an ODE
solver, we utilize a classical explicit midpoint time-integration
method with a time-step of $2^{-12} \approx 2.4 \times 10^{-4}$.  In the full
system, we evaluate right-hand side the via a fast method
from~\cite{matveev2015a-fast8543937}.

\begin{figure}
  \centering
  \includegraphics{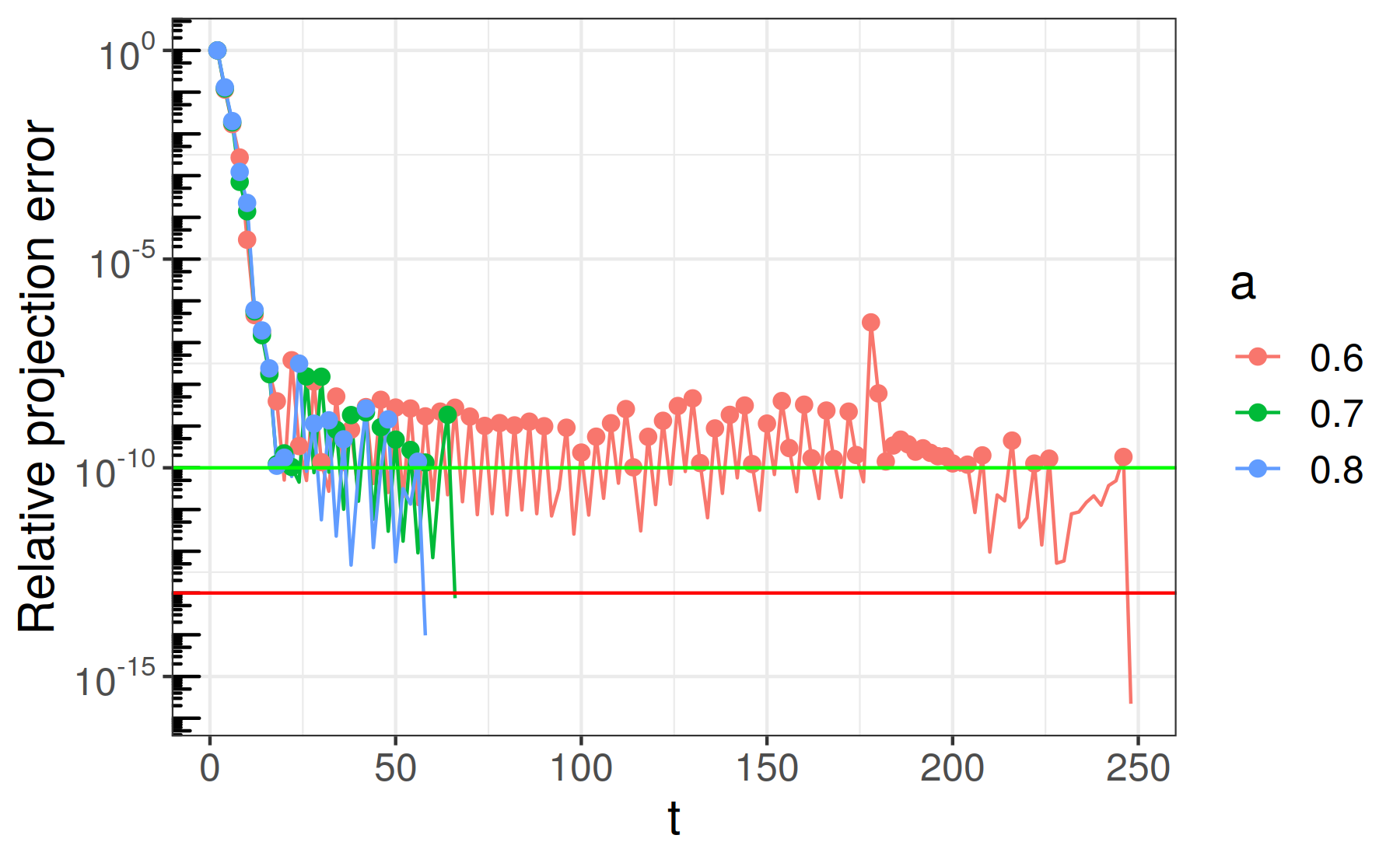}

  \caption[Projection error]{The dependency of the new basis
    projection error on time.  On the vertical axis are the
    $\left\lVert (I - V_{k-1} V_{k-1}^T) \hat{V}_k \right\rVert_2$
    values; dots denote the moments when the basis is expanded.  The
    horizontal green line shows the $\varepsilon'$ value, the red one~---
    $\varepsilon$.}\label{fig:proj}
\end{figure}

Figure~\ref{fig:proj} demonstrates the inner working process of the
algorithm during the basis construction for kernel parameter $a$
equal to $0.6$, $0.7$ or $0.8$.  We can see that the projection error
decreases rapidly at the onset, then oscillates a bit around
$\varepsilon'$, and eventually crosses the $\varepsilon$ boundary.
This, incidentally, highlights another role of the $\varepsilon'$
parameter~--- it effectively prevents the algorithm from
over-approximating an initial segment of the solution; we show
below the reason why this is important. 

\begin{figure}
  \centering
  \includegraphics{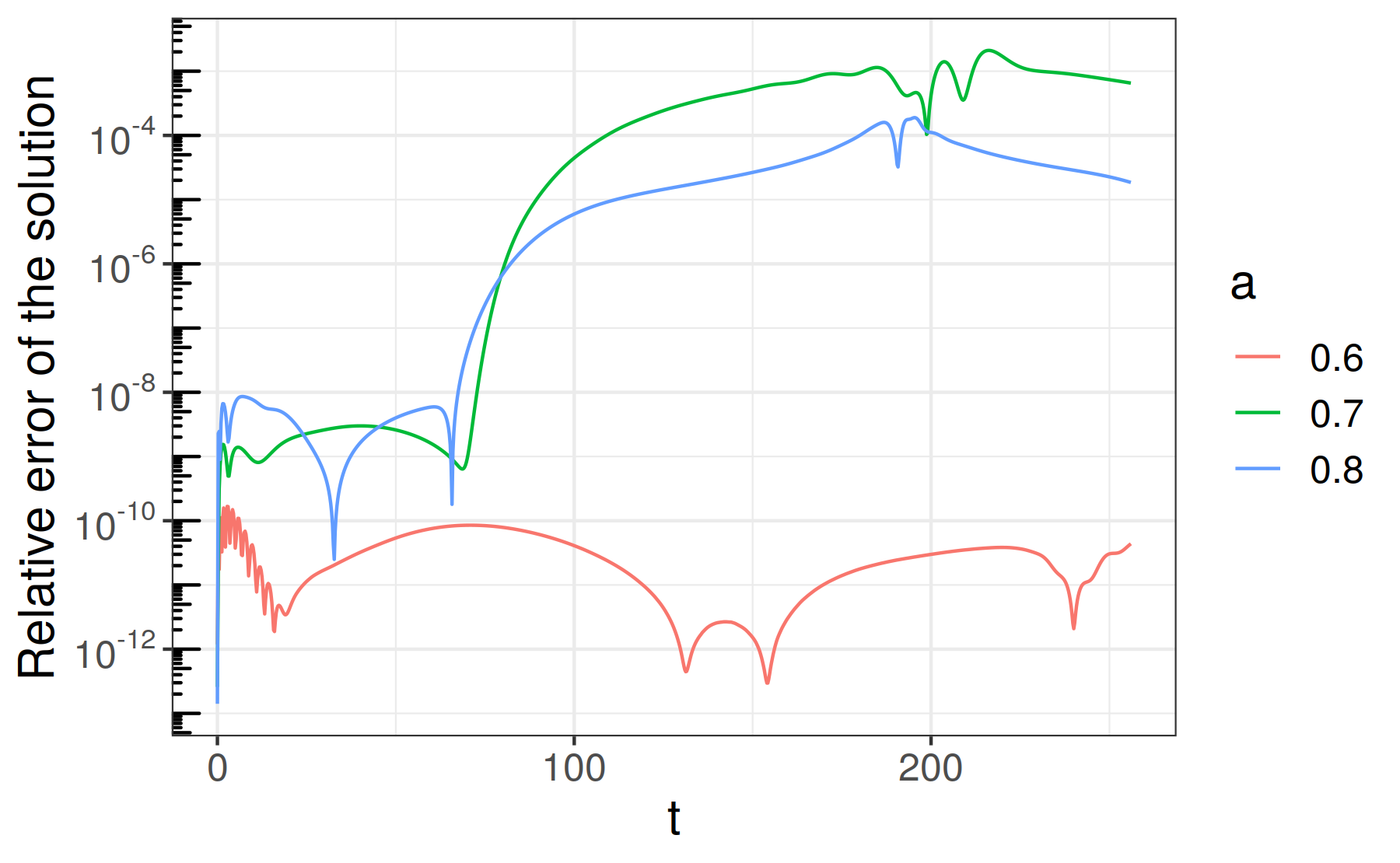}

  \caption[Solution error]{The dependency of the reduced solution
    relative error in Euclidean norm on time for $N = 32768$.  On the vertical axis are
    the values of ${\left\lVert n(t) - \tilde{n}(t) \right\rVert}_2 /
    {\left\lVert n(t) \right\rVert}_2$, where $\tilde{n}(t) = V
    \tilde{x}(t)$.}\label{fig:sol}
\end{figure}

Figure~\ref{fig:sol} demonstrates an error for a recomputed
solution of a reduced system~(\ref{eq:smol_reduced}) for the
time segment $t \in [0,256]$, as compared to a solution of a full
system~(\ref{eq:smoluchowski}).  As can be seen from the figure, the
error remains small on an interval where the basis was originally
constructed, and sharply increases at its end (where the algorithm
effectively switches from interpolation to extrapolation)~--- but,
crucially, it still remains bounded around acceptable level $10^{-3}$.

\begin{table}
  \centering
  \begin{tabular}{rrrr}
    \toprule
    $a$ & $t_\mathrm{full}$, sec & $t_\mathrm{red}$, sec & Basis size \\
    \midrule
    $0.6$ & $5.3 \times 10^3$ & $10^4$ & 216 \\
    $0.7$ & $5.3 \times 10^3$ & $290$ & 99 \\
    $0.8$ & $5.3 \times 10^3$ & $61$ & 86 \\
    \bottomrule
  \end{tabular}
  \caption{Solution time for $N = 32768$.  $t_\mathrm{full}$ is time to solve the full
    system~(\ref{eq:smoluchowski}), $t_\mathrm{red}$---time to
    solve~(\ref{eq:smol_reduced}).}\label{tab:timings}
\end{table}

For the same set of simulations, Table~\ref{tab:timings} provides the CPU
time required to solve full and reduced systems, as well as the
eventual size of the basis.  As the table clearly demonstrates, the
use of the reduction is not always beneficial, especially if fast
algorithms for evaluation of the full operator are available;
specifically, in the case $a = 0.6$, with the basis size of 216, the
reduced system is actually more expensive to solve than the full problem.

Since this observation, together with Figure~\ref{fig:sol},
strongly hint at a trade-off between performance and precision, one
might be tempted to tweak the $\varepsilon$ parameter to manage it.
Unfortunately, as our second set of experiments demonstrates, this is
not always straight-forward in practice, which we find rather 
surprising.

\begin{table}
  \centering
  \begin{tabular}{rrrrr}
    \toprule
    $a$ & $\varepsilon$ & Reduced solution error & Time span used for basis & Basis size \\
    \midrule
    $0.7$ & $10^{-10}$ & $3.5 \times 10^{-3}$ & $[0,100]$ & 68 \\
    $0.7$ & $10^{-11}$ & $2.3 \times 10^{-1}$ & $[0,28]$ & 52 \\
    $0.7$ & $10^{-12}$ & $2.8 \times 10^{-5}$ & $[0,128]$ & 101 \\
    $0.7$ & $10^{-13}$ & $2.3 \times 10^{-2}$ & $[0,68]$ & 102 \\
    $0.7$ & $5 \times 10^{-14}$ & $4.2 \times 10^{-3}$ & $[0,76]$ & 112 \\
    $0.7$ & $10^{-14}$ & $4.6 \times 10^{-10}$ & $[0,256]$ & 234 \\
    \bottomrule
  \end{tabular}
  \caption{Solution error as a function of $\varepsilon$ for
    $N = 65536$.  Reduced solution error is taken as a maximum
    relative error in 2-norm over the time span of
    $[0,256]$.}\label{tab:solerr_eps}
\end{table}

These experiments are performed with $N = 65 536$ (we find that the effect is
more visible at this dimensionality), with $\delta = \varepsilon$ and
$\varepsilon' = 10^3 \cdot \varepsilon$. The results are available in
Table~\ref{tab:solerr_eps}; all the other parameters are the same as
in the case above.  As can be readily observed, the resulting error
does not depend monotonically on $\varepsilon$, and seems to depend
more on the actual time span which was used to construct the basis.
Note that in the very last row, corresponding to $\varepsilon =
10^{-14}$, the algorithm has simply used up the entire time span under
evaluation for basis construction, and therefore the error reflects
the \enquote{interpolation} mode, as seen in Figure~\ref{fig:sol}.

Finally, to test the scalability of our approach, we have tested the
algorithm with a larger system with $N = 131072$ and
$a \in \{0.6, 0.7\}$.

\begin{figure}
  \centering
  \includegraphics{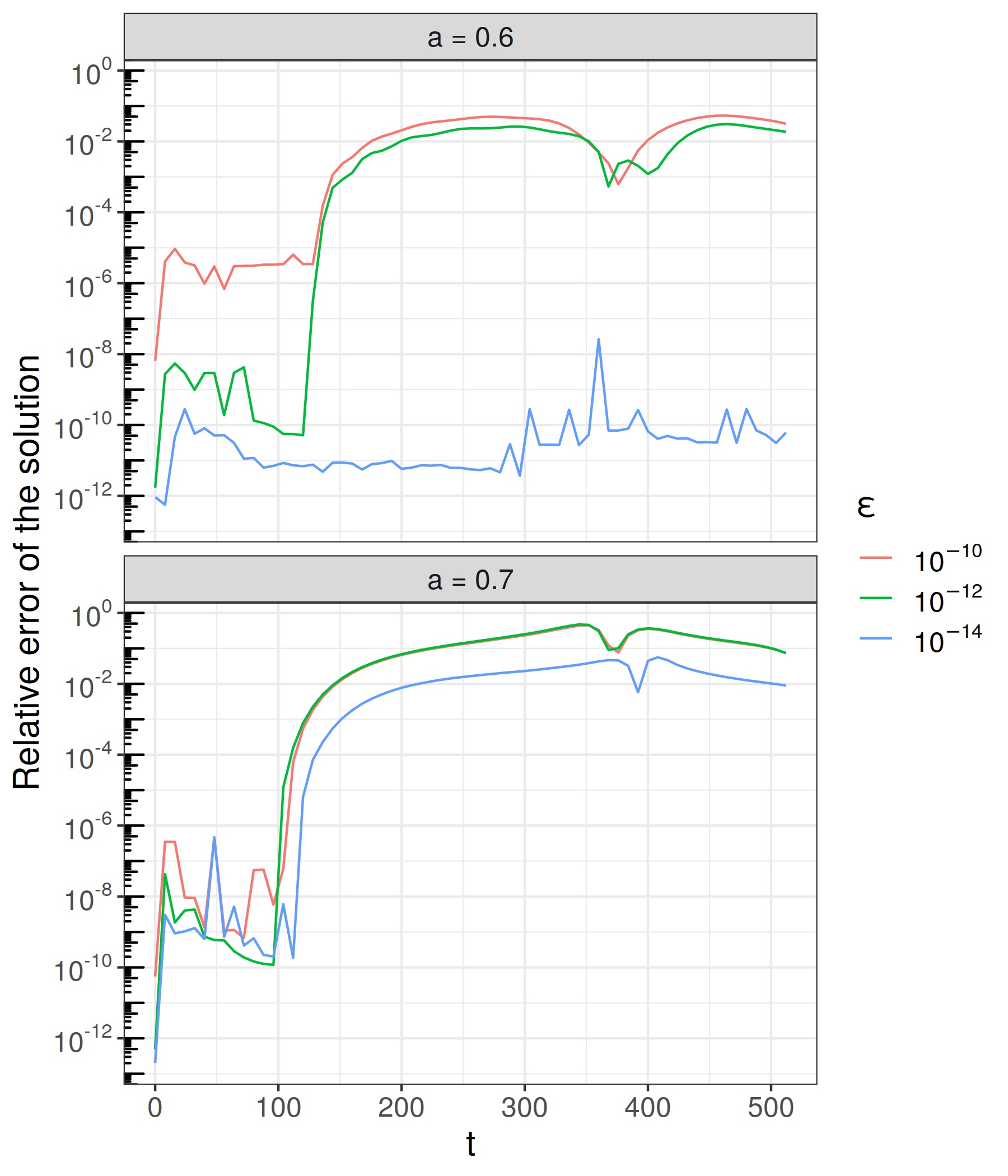}
  \caption{Solution error for $N = 131072$, $a = 0.7$ and $a = 0.6$}\label{fig:sol_large_err}
\end{figure}

\begin{table}
  \centering
  \begin{tabular}{rrrrrr}
    \toprule
    $a$ & $\varepsilon$ & $t_\mathrm{full}$, sec & $t_\mathrm{red}$, sec & Basis size & Time span used for basis \\
    \midrule
    $0.6$ & $10^{-10}$ & $8.8 \times 10^4$ & $233.18$ & 84 & $[0,128]$ \\
    $0.6$ & $10^{-12}$ & $8.8 \times 10^4$ & $3.3 \times 10^3$ & 115 & $[0, 124]$ \\
    $0.6$
        & $10^{-14}$ & $4.5 \times 10^4$ & $1.2 \times 10^5$ & 459 & $[0,508]$ \\
    \midrule
    $0.7$ & $10^{-10}$ & $8.8 \times 10^4$ & $144$ & 70 & $[0,106]$ \\
    $0.7$ & $10^{-12}$ & $8.8 \times 10^4$ & $1.3 \times 10^3$ & 93 & $[0, 94]$ \\
    $0.7$ & $10^{-14}$ & $8.8 \times 10^4$ & $6.2 \times 10^3$ & 151 & $[0,104]$ \\
    \bottomrule
  \end{tabular}
  \caption{Solution time and basis information for $N = 131072$; the test with $a = 0.6$ and $\varepsilon = 10^{-14}$ was run on a different hardware than the rest.}\label{tab:timings_large}
\end{table}

Figure~\ref{fig:sol_large_err} demonstrates the relative error of the
reconstructed solution for $\varepsilon = 10^{-10}$,
$\varepsilon = 10^{-12}$ and $\varepsilon = 10^{-14}$, on a longer
time segment of $t \in [0,512]$.  The same sharp
\enquote{interpolation---extrapolation} transition is visible here;
and results from Table~\ref{tab:timings_large} confirm that, in both
cases, the \enquote{interpolation} region is in the neighbourhood of
the transition visible on the graph, except for $a = 0.6$,
$\varepsilon = 10^{-14}$, where almost the entire time span was used
for basis construction, and thus there is no visible transition at
all.

\begin{figure}
  \centering
  \includegraphics{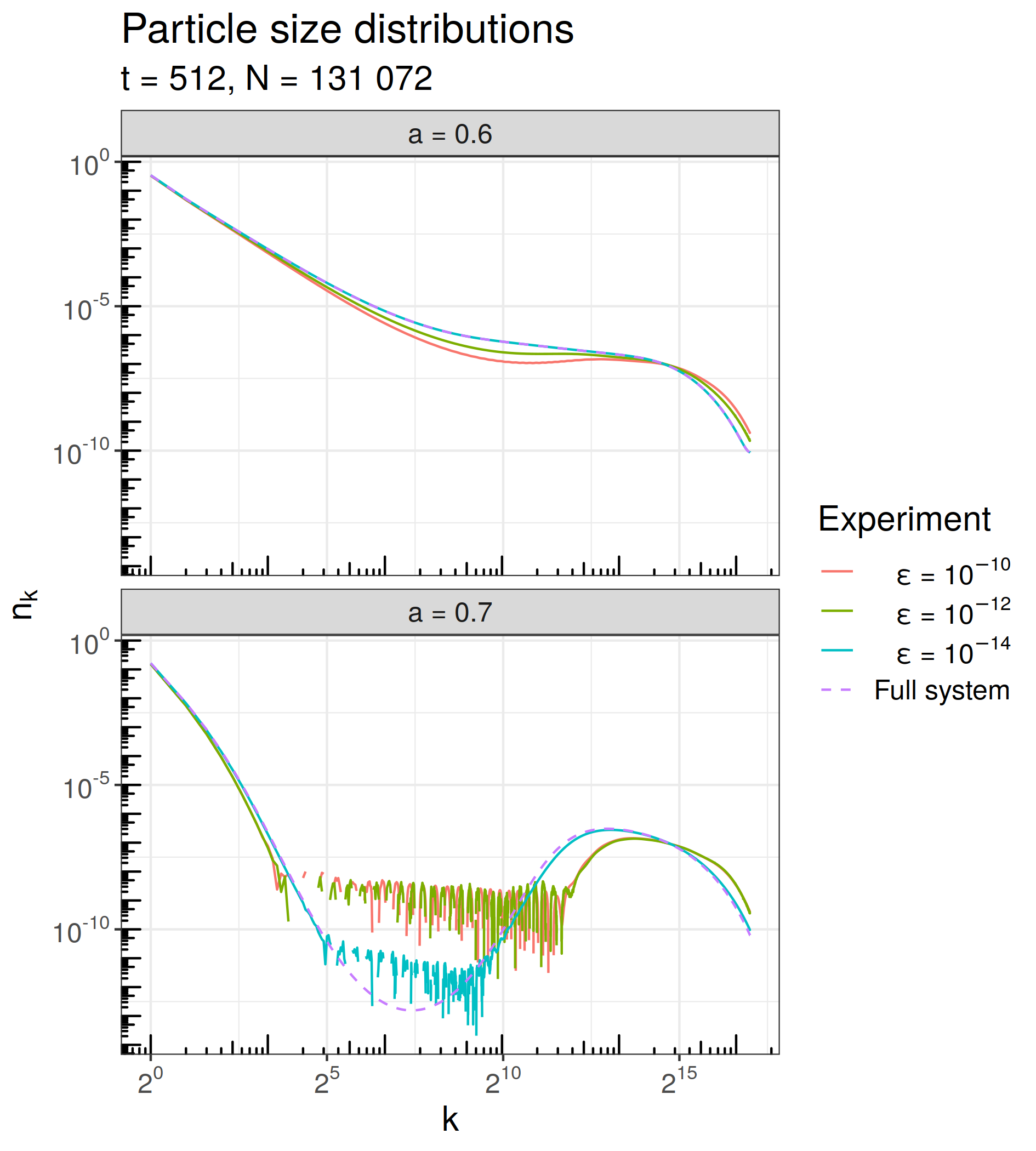}
  \caption{Full particle size distribution (purple line) and reduced
    solutions (red, green and blue) at $t = 512$ for $N = 131072$, $a = 0.7$
    and $a = 0.6$. The solution with $\varepsilon = 10^{-14}$ is close
    to the full solution, diverging only for the smallest
    concentration values}\label{fig:sol_large}
\end{figure}

Finally, Figure~\ref{fig:sol_large} demonstrates full and reduced
solutions at the far end of the simulation time-interval.  The
solution with $a = 0.6$, $\varepsilon = 10^{-14}$ is indistinguishable
from the precise one, as it is effectively computed in
\enquote{interpolation} mode, but even the solution with $a = 0.7$,
$\varepsilon = 10^{-14}$ is fairly quantitatively close to the full
solution, diverging only for the smallest values concentration of
$n_k$~--- which, surprisingly, do not affect significantly the
solution to the either side of their mass range.

Less precise solutions with $\varepsilon = 10^{-10}$ and
$\varepsilon = 10^{-12}$ do not deliver a good quantitative fit (as
can already be seen in Figure~\ref{fig:sol_large_err}) but, nevertheless,
reproduce the qualitative shape of the solution well, despite
being significantly cheaper to compute in terms of CPU-time.

\section{Conclusion}\label{sec:conclusion}

We have suggested an application of the popular and well-established method of
model reduction to the problem of a system of Smoluchowski ODEs,
including a candidate method for construction of a reduced basis
from an automatically selected prefix of the modelled time span.
In our numerical experiments, we demonstrate the existence of such
a low-dimensional basis, noticeable speed-ups  of computations, and
 reasonable approximation to the full solution by the reduced model. 

At the same time, we also demonstrate problematic sides of the chosen approach~--- 
the  precision control of the reduced solution seems to be not straight-forward 
at all due to the nonlinearity of both method and model. Hence, control of accuracy requires more theoretical 
analysis. In light of these shortcomings, we find our concept 
very promising for future development in more complicated applied cases and also consider it as 
fruitful directions of further research.
 
\section{Acknowledgements}

We are grateful to Nikolai Zamarashkin for comprehensive discussions
during preparation of this work. The work was supported by the Russian
Science Foundation, grant~\mbox{19--11--00338}.

\bibliography{reduction.en}

\end{document}